\documentclass{article}
\usepackage[utf8]{inputenc}
\usepackage[margin=1in]{geometry}
\usepackage{mathtools}
\usepackage{amsmath}
\usepackage{amssymb}
\usepackage{amsthm}
\usepackage{enumitem}
\usepackage{hyperref}

\newtheorem{lemma}{Lemma}
\newtheorem{theorem}{Theorem}

\newcommand{\C}{\mathbb{C}}

\DeclarePairedDelimiter{\abs}{\lvert}{\rvert}
\DeclareMathOperator{\sgn}{sgn}

\title{An Elementary Proof of the Signature of Satellite Knots}
\author{Daniel Carter}
\date{}

\begin{document}

\maketitle

\abstract{We present a proof of Litherland's formula for the Tristram-Levine signature of a satellite knot in terms of its constituents. Litherland's original proof used more advanced algebraic techniques, while ours uses only linear algebra and some basic results in knot theory.}

\section{Background}

First, recall the definition of a satellite knot and the Tristram-Levine signature. A (nontrivial) \textit{satellite knot} is obtained as follows:
\begin{enumerate}
    \item Embed a knot $K$ in the solid torus $T=S^1\times D^2$. For nontrivial satellite knots we require that there is no simply-connected subset of the solid torus containing the knot and $K$ is not isotopic (in the solid torus) to the central $S^1$ of the solid torus.
    \item Take the image of $K$ under a homeomorphism taking $T$ to the (solid) tubular neighborhood of another knot $J$. We require that this homeomorphism is ``untwisted,'' that is, linking numbers between any two closed curves in $T$ are preserved in the image.
\end{enumerate}

Here, $K$ together with the embedding into $T$ is called the \textit{pattern} and $J$ is called the \textit{companion knot}. A special class of satellites are those with $K=$ the $(p,q)$ torus knot embedded in the standard way onto the surface of the torus. This is called the \textit{$(p,q)$ cable} of $J$ (see Chapter 1 in \cite{lickorish}).

The \textit{Tristram-Levine signature} is a knot invariant defined as the signature of the matrix $(1-\omega)M-(1-\overline{\omega})M^T$, where $M$ is any Seifert matrix of the knot and $\omega$ is a complex number with $\abs{\omega}=1$. The fact that this signature is the same for all Seifert surfaces of a knot follows by considering the effect from performing surgery along an arc to transform one Seifert surface to another; see Theorem 8.9 in \cite{lickorish}. The signature of a knot $K$ is denoted by $\sigma_\omega(K)$.

We prove the following useful formula relating the Tristram-Levine signature of a satellite knot to the signatures of the constituent knots:

\begin{theorem}\label{signaturesatellite}
If $K'$ is a satellite of $J$ by $K$ and $n$ is the winding number of the embedding of $K$ in the solid torus, then
\[ \sigma_{\omega}(K')=\sigma_{\omega}(K)+\sigma_{\omega^n}(J). \]
\end{theorem}

This formula was proven by Litherland in 1979 (\cite{litherland}) in order to study algebraic knots, which are a subset of the set of cables of cables of ... of torus knots. This was spurred by Rudolph's question about the independence of algebraic knots in the concordance group \cite{rudolph}. Litherland's proof of the formula for the signature of a satellite knot uses algebraic techniques. We provide a proof which uses only linear algebra and some basic results of knot theory.

\section{Lemmas}
\label{lemmas}
We will need the following four lemmas. First, a fact of linear algebra:

\begin{lemma}\label{matrixlem}
Suppose $M$ is a Hermitian matrix and $M'$ is obtained from $M$ by one of the following operations:
\begin{itemize}
    \item Add $z$ times row $r_i$ to $r_j$, then $\overline{z}$ times column $c_i$ to $c_j$ for some $z\in\C$.
    \item Replace row $r_i$ with $zr_i$, then column $c_i$ with $\overline{z}c_i$ for some $z\in\C\setminus\{0\}$.
\end{itemize}
Then $M'$ is also Hermitian, and if $\det(M)\ne 0$ then $\det(M')\ne 0$ and $\sgn(M)=\sgn(M')$.
\begin{proof}
This is a special case of Sylvester's law of inertia for complex matrices; note that both of the operations above are congruences.
\end{proof}
\end{lemma}

Then we will need the following three facts from knot theory. All of these may be found in \cite{lickorish}.

\begin{lemma}[Part of the proof of Theorem 6.10(ii) in \cite{lickorish}]\label{blockdiagonal}
For an appropriate choice of generators of the Seifert surface homology, if $A$ is the Seifert matrix associated to those generators, then
\[ A-A^T=\begin{bmatrix}
0 & 1 & 0 & 0 & \dots & 0 & 0 \\
-1 & 0 & 0 & 0 & \dots & 0 & 0 \\
0 & 0 & 0 & 1 & \dots & 0 & 0 \\
0 & 0 & -1 & 0 & \dots & 0 & 0 \\
\vdots & \vdots & \vdots  & \vdots & \ddots & \vdots & \vdots \\
0 & 0 & 0 & 0 & \dots & 0 & 1 \\
0 & 0 & 0 & 0 & \dots & -1 & 0
\end{bmatrix} \]
where the number of $\begin{bmatrix}0&1\\ -1&0\end{bmatrix}$ blocks on the diagonal is equal to the genus of the Seifert surface.
\end{lemma}

\begin{lemma}[Part of the proof of Theorem 6.15 in \cite{lickorish}]\label{seifertmatrix}
If $K'$ is a satellite of $J$ by $K$ and $n$ is the winding number of the embedding of $K$ in the solid torus, then a Seifert matrix $A$ for $K'$ is given by the block matrix
\[ A=\begin{bmatrix}
M & 0 & 0 & \cdots & 0 \\
0 & N & N & \cdots & N \\
0 & N^T & N & \cdots & N \\
\vdots & \vdots & \vdots & \ddots & \vdots \\
0 & N^T & N^T & \cdots & N
\end{bmatrix} \]
where $M$ is any Seifert matrix for $K$ and $N$ any Seifert matrix for $J$, and there are $n\times n$ copies of $N$ and $N^T$.
\end{lemma}

\begin{lemma}[Theorem 6.15 in \cite{lickorish}]\label{alexanderpoly}
If $K'$ is a satellite of $J$ by $K$ and $n$ is the winding number of the embedding of $K$ in the solid torus, then
\[ \Delta_{K'}(t)\doteq \Delta_K(t)\Delta_J(t^n) \]
where $\Delta$ is the Alexander polynomial and $\doteq$ means equal up to a factor of $t^k$.
\end{lemma}

\section{Proof of Theorem~\ref{signaturesatellite}}
\label{signature}

The proof of Lemma~\ref{alexanderpoly} essentially amounts to three steps after pulling out the $\det(tM-M^T)$ term from the matrix obtained by Lemma~\ref{seifertmatrix}:
\begin{enumerate}
    \item Obtain $t^nN-N^T$ in the matrix by row operations.
    \item Zero out the other blocks in the same row as the $t^nN-N^T$ by column operations.
    \item Pull out the $\det(t^nN-N^T)$ term and compute that the determinant of what remains is a unit.
\end{enumerate}

Our proof is similar, though it is more difficult because the signature is only preserved under congruences (as opposed to the determinant, which changes predictably with arbitrary row/column operations). In particular, we do the following after pulling out a $\sgn((1-\omega)M+(1-\overline{\omega})M^T)$ term:
\begin{enumerate}
    \item Obtain (a multiple of) $(1-\omega^n)N+(1-\overline{\omega}^n)N^T$ in the matrix by congruence.
    \item Zero out the other blocks in the same row and column as $(1-\omega^n)N+(1-\overline{\omega}^n)N^T$ by congruence.
    \item Pull out the $\sgn((1-\omega^n)N+(1-\overline{\omega}^n)N^T)$ term and compute that the signature of the remaining matrix is zero.
\end{enumerate}

In the first and third steps, there are some additional subtleties not found in the proof of Lemma~\ref{alexanderpoly}. For the first step, there are some exceptional $\omega$ which make the determinant vanish; these are dealt with by applying Lemma~\ref{alexanderpoly}. For the third step, we involve an additional congruence to get the matrix into a form where the signature may be readily calculated using Lemma~\ref{blockdiagonal}.

Now here is the proof of Theorem~\ref{signaturesatellite}:

\begin{proof}
Let $A$ be the Seifert matrix for $K'$ described in Lemma~\ref{seifertmatrix}. Now consider
\begin{align*}
    (1-\omega)A&+(1-\overline{\omega})A^T = ((1-\omega)M+(1-\overline{\omega})M^T) \\
    &\oplus\begin{bmatrix}
    (1-\omega)N+(1-\overline{\omega})N^T & (1-\omega)N+(1-\overline{\omega})N & \cdots & (1-\omega)N+(1-\overline{\omega})N \\
    (1-\omega)N^T+(1-\overline{\omega})N^T & (1-\omega)N+(1-\overline{\omega})N^T & \cdots & (1-\omega)N+(1-\overline{\omega})N \\
    \vdots & \vdots & \ddots & \vdots \\
    (1-\omega)N^T+(1-\overline{\omega})N^T & (1-\omega)N^T+(1-\overline{\omega})N^T & \cdots & (1-\omega)N+(1-\overline{\omega})N^T
    \end{bmatrix},
\end{align*}
with $\oplus$ the direct sum of matrices. The signature of this is the signature of $(1-\omega)M+(1-\overline{\omega})M^T$ plus the signature of the second block matrix, which we will call $B$. We will perform some carefully chosen row/column operations on $B$ so that the signature does not change, guaranteed by Lemma~\ref{matrixlem}.

Let $X$ be the (block) row matrix given by $\sum_{k=1}^n (\omega^{1-k}+\dots+\omega^{n-k})\times(\text{block row $k$ of $B$})$. Let $Y$ be the (block) column matrix given by $\sum_{k=1}^n (\overline{\omega}^{1-k}+\dots+\overline{\omega}^{n-k})\times(\text{block column $i$ of $B$})$. Observe by a simple telescoping argument that all the blocks of $X$ and $Y$ are just $(1-\omega^n)N+(1-\overline{\omega}^n)N^T$. Replace the first row and column of $B$ with $X$ and $Y$, except the top-left block, which is twice $(1-\omega^n)N+(1-\overline{\omega}^n)N^T$; this corresponds to doing all the row operations necessary to make the first row $X$ and all the column operations to make the first column $Y$. Order these operations so you first replace block row 1 with
\[ 1+\omega+\dots+\omega^{n-1} \]
times itself and block column 1 with
\[ 1+\overline{\omega}+\dots+\overline{\omega}^{n-1} \]
times itself. This corresponds to a matrix operation of the second type from Lemma~\ref{matrixlem} as long as $1+\omega+\dots+\omega^{n-1}\ne 0$, in which case the signature is unchanged. This value is zero exactly when $\omega$ is an $n$th root of unity other than 1. But if $\omega$ is an $n$th root of unity, note by factoring $(1-\omega)A+(1-\overline{\omega})A^T=(1-\omega)(A-\overline{\omega}A^T)$ that the value of $\sigma$ as a function of $\omega$ can only change at zeros of the Alexander polynomial, and by Lemma~\ref{alexanderpoly}, $\Delta_{K'}(t)\doteq\Delta_{K}(t)\Delta_{J}(t^n)$. Set $t=\omega$ and note that $\Delta_{J}(\omega^n)=\Delta_{J}(1)\ne 0$ so $\Delta_{J}(t^n)$ does not have a zero at any $n$th root of unity, which is sufficient for this case.

The other row/column operations are the first type from Lemma~\ref{matrixlem}, so they also do not change the signature.

Our matrix now looks like
\[ \begin{bmatrix}
2((1-\omega^n)N+(1-\overline{\omega}^n)N^T) & (1-\omega^n)N+(1-\overline{\omega}^n)N^T & \cdots & (1-\omega^n)N+(1-\overline{\omega}^n)N^T \\
(1-\omega^n)N+(1-\overline{\omega}^n)N^T & (1-\omega)N+(1-\overline{\omega})N^T & \cdots & (1-\omega)N+(1-\overline{\omega})N \\
\vdots & \vdots & \ddots & \vdots \\
(1-\omega^n)N+(1-\overline{\omega}^n)N^T & (1-\omega)N^T+(1-\overline{\omega})N^T & \cdots & (1-\omega)N+(1-\overline{\omega})N^T
\end{bmatrix}.
\]

Subtract half of row 1 and half of column 1 from the rest of the matrix. The result is
\[ (2((1-\omega^n)N+(1-\overline{\omega}^n)N^T))\oplus\begin{bmatrix}
D & U & U & \cdots & U \\
L & D & U & \cdots & U \\
L & L & D & \cdots & U \\
\vdots & \vdots & \vdots & \ddots & \vdots \\
L & L & L & \cdots & D
\end{bmatrix}
\]
with
\begin{align*}
    D &= (\omega^n-\omega)N+(\overline{\omega}^n-\overline{\omega})N^T \\
    U &= (\omega^n-\omega-\overline{\omega}+1)N + (\overline{\omega}^n-1)N^T \\
    L &= (\omega^n-1)N + (\overline{\omega}^n-\overline{\omega}-\omega+1)N^T.
\end{align*}
The signature is unchanged, and the latter matrix in the direct sum above is now $(n-1)\times(n-1)$ blocks; e.g. the first block row is one $D$ and $(n-2)$ $U$'s.

So the signature of the matrix in question is the signature of $2((1-\omega^n)N+(1-\overline{\omega}^n)N^T)$ plus the signature of the rest, which we will call $C$. Subtract $\frac{1}{2(n-2)}$ times each block row and column of $C$ from every other block row and column. This does not change the signature. The new diagonal blocks are equal to
\[ \frac{1}{2}(-\omega+\overline{\omega})(N-N^T). \]
The blocks immediately above the diagonal are
\[ \left(-\frac{n-3}{2(n-2)}\omega-\frac{n-1}{2(n-2)}\overline{\omega}+\frac{n-2}{n-2}\right)(N-N^T). \]
The blocks above that are
\[ \left(-\frac{n-4}{2(n-2)}\omega-\frac{n-2}{2(n-2)}\overline{\omega}+\frac{n-3}{n-2}\right)(N-N^T) \]
since they incorporate one more $U$ term and one less $L$ term. This pattern continues, so that a block $d$ above the diagonal is now
\[ \left(-\frac{n-2-d}{2(n-2)}\omega-\frac{n-d}{2(n-2)}\overline{\omega}+\frac{n-1-d}{n-2}\right)(N-N^T). \]
The matrix is still Hermitian, so a block $d$ below the diagonal is
\[ \left(\frac{n-2-d}{2(n-2)}\overline{\omega}+\frac{n-d}{2(n-2)}\omega-\frac{n-1-d}{n-2}\right)(N-N^T). \]
This means the new matrix is equal to the Kronecker product
\[ \begin{bmatrix}
\frac{1}{2}(-\omega+\overline{\omega}) & -\frac{n-3}{2(n-2)}\omega-\frac{n-1}{2(n-2)}\overline{\omega}+\frac{n-2}{n-2} & \cdots \\
\frac{n-3}{2(n-2)}\overline{\omega}+\frac{n-1}{2(n-2)}\omega-\frac{n-2}{n-2} & \frac{1}{2}(-\omega+\overline{\omega}) & \cdots \\
\vdots & \vdots & \ddots
\end{bmatrix}\otimes (N-N^T). \]
The eigenvalues of this are equal to the pairwise products of eigenvalues of the two matrices. Notice that the first matrix is skew-Hermitian, so its eigenvalues are all pure imaginary. By Lemma~\ref{blockdiagonal}, the latter matrix is of the form
\[ \begin{bmatrix}
0 & 1 & 0 & 0 & \dots & 0 & 0 \\
-1 & 0 & 0 & 0 & \dots & 0 & 0 \\
0 & 0 & 0 & 1 & \dots & 0 & 0 \\
0 & 0 & -1 & 0 & \dots & 0 & 0 \\
\vdots & \vdots & \vdots  & \vdots & \ddots & \vdots & \vdots \\
0 & 0 & 0 & 0 & \dots & 0 & 1 \\
0 & 0 & 0 & 0 & \dots & -1 & 0
\end{bmatrix} \]
for appropriate choice of generators of the Seifert surface homology. This has eigenvalues $i$ and $-i$ with equal multiplicities. Hence the eigenvalues of the Kronecker product come in $+$/$-$ pairs with equal multiplicity, so the signature of $C$ is 0.

We have shown
\begin{align*}
    \sgn((1-\omega)A+(1-\overline{\omega})A^T) &= \sgn((1-\omega)M+(1-\overline{\omega})M^T) + \sgn(2((1-\omega^n)N+(1-\overline{\omega}^n)N^T)) \\
    &= \sgn((1-\omega)M+(1-\overline{\omega})M^T) + \sgn((1-\omega^n)N+(1-\overline{\omega}^n)N^T) \\
    \sigma_{\omega}(K') &= \sigma_{\omega}(K) + \sigma_{\omega^n}(J).
\end{align*}
\end{proof}

\section{Closing Remarks}

One might compare this proof to \cite{shinohara}, in which the special case $\omega=-1$ is proven also using essentially only linear algebra. It appears, however, that Shinohara's proof does not generalize; in particular, the structure of the matrix obtained after the congruence fundamentally depends on the parity of $n$, so Shinohara obtains specifically the formula
\[ \sigma_{-1}(K')=\begin{cases}
\sigma_{-1}(K) & \text{if $n$ is even}, \\
\sigma_{-1}(K)+\sigma_{-1}(J) & \text{if $n$ is odd}
\end{cases} \]
as a result. For other $\omega$, the signature formula depends on more than just the parity of $n$, so a different matrix congruence must be used.

Litherland also remarks that the result holds if $K$ is a link as well, with the exception when $\omega$ is an $n$th root of unity, where the formula may be off by up to $\pm2(m-1)$ where $m$ is the number of components of $K$. This is clear in our proof as well; the exception occurs when $1+\omega+\dots+\omega^{n-1}$ is zero (that is, $\omega$ is an $n$th root of unity other than $1$) because the use of Theorem 6.15 from \cite{lickorish} is only justified when $K$ is a knot, not a link.

\section{Acknowledgements}

Thanks to Ian Zemke for identifying some minor errors in the proof of Theorem~\ref{signaturesatellite} and Ollie Thakar for feedback and encouragement.

The author is extremely grateful to Jae Choon Cha for bringing his paper \cite{cha} to the author's attention after the initial arXiv submission of this paper, and for his kind words. See the addendum below for a discussion of \cite{cha} as it relates to this work.

\section{Addendum}

In Lemma 2.2 in \cite{cha}, it is shown that what we call matrix $B$ is congruent to (in our notation)
\[ \frac{\omega^n N - N^T}{\omega^n-1}\oplus\bigoplus_{k=1}^n \frac{\omega^{k+1}-1}{(\omega-1)(\omega^k-1)}(N - N^T). \]
The first term here is $\frac{1}{2-\omega-\overline{\omega}}$ (a real number) times $(1-\omega^n)N+(1-\omega^n)N^T$, so they have same signature. The signature of all the other summands is 0. Cha and Ko find this congruence iteratively, one term at a time. Contrast this with our method, where we do our whole congruence in a few large steps.

Note that the large direct sum above is equal to the Kronecker product
\[ \Delta\otimes (N-N^T) = \begin{bmatrix}
\frac{\omega^2-1}{(\omega-1)(\omega-1)} & 0 & 0 & \cdots \\
0 & \frac{\omega^3-1}{(\omega-1)(\omega^2-1)} & 0 & \cdots \\
0 & 0 & \frac{\omega^4-1}{(\omega-1)(\omega^3-1)} & \cdots \\
\vdots & \vdots & \vdots & \ddots
\end{bmatrix}\otimes (N-N^T) \]
where $\Delta$ is a pure imaginary diagonal matrix. Our matrix $C$ was itself congruent to a Kronecker product of a skew-Hermitian matrix $S$ and $(N-N^T)$. Hence diagonalizing $iS$ (a Hermitian matrix) and scaling the diagonal entries appropriately should result in $i\Delta$ (a real diagonal matrix). If you do this diagonalization process entry-by-entry and compose that process with our congruence, you should get exactly Cha and Ko's method.

The result in \cite{cha} is actually more general in two main ways:
\begin{itemize}
    \item Let $\varepsilon\in\{-1,1\}$. Replace $N^T$ by $\varepsilon N^T$ everywhere. The convention that the signature of a skew-Hermitian matrix is equal to the signature of $i$ times the matrix. The resulting signature formula includes a term involving $\sgn(N-\varepsilon N^T)$; in the $\varepsilon=1$ case (our case) this signature is zero so this term does not appear.
    \item Replace some $n-u$ diagonal blocks of $B$ with $(1-\omega)\varepsilon N^T+(1-\overline{\omega})N$. By rearranging the rows and columns, these can be taken to be the last $n-u$ blocks, so blocks $1$ through $u$ are unchanged. The resulting signature formula replaces $n$ with $2u-n$.
\end{itemize}
The iterative method in \cite{cha} lends itself well to these generalizations. Our proof can also handle these generalizations to an extent. In particular:
\begin{itemize}
    \item If $\varepsilon=-1$ then exactly the same congruences work, just everywhere replacing $N^T$ with $-N^T$. The resulting formula is $\sgn(B)=\sgn((1-\omega^n)N-(1-\overline{\omega}^n)N^T)+\sgn(S)\sgn(N+N^T)$. It is easy to verify that $\sgn(A\otimes B)=\pm \sgn(A)\sgn(B)$ where $A$ and $B$ are (skew-)Hermitian, taking the $+$ sign if at least one is Hermitian and the $-$ sign if both are skew-Hermitian. Unfortunately, it does not seem there is an easy way to compute $\sgn(S)$ other than by diagonalization. The answer, in \cite{cha}, is $\sgn(S)=n+1-2\lceil{\frac{nx}{2\pi}}\rceil$ where $\omega=e^{ix}$ with $0\le x\le 2\pi$, unless $\frac{nx}{2\pi}$ is an integer in which case it is this expression minus 1. At $x=0$, $\sgn(S)=0$. This may be seen fairly easily from $\Delta$, but given only $S$ it is somewhat surprising the signature has such a simple expression.
    \item If the last $n-u$ diagonal blocks are replaced with $\varepsilon((1-\omega)N^T+(1-\overline{\omega})N)$, then replace $X$ in the first step with
    \[ \sum_{k=1}^u(\omega^{1-k}+\dots+\omega^{2u-n-k})\times(\text{block row $k$ of $B$}) - \sum_{k=u+1}^n(\omega^{k-2u}+\dots+\omega^{k-n-1})\times(\text{block row $k$ of $B$}) \]
    and modify $Y$ similarly. After this, perform the same step of subtracting half of row 1 and half of column 1 from the rest. The result is now
    \[ (2((1-\omega^{2u-n})N+(1-\overline{\omega}^{2u-n})N^T))\oplus\begin{bmatrix}
    D_2 & U & U & \cdots & U \\
    L & D_3 & U & \cdots & U \\
    L & L & D_4 & \cdots & U \\
    \vdots & \vdots & \vdots & \ddots & \vdots \\
    L & L & L & \cdots & D_n
    \end{bmatrix}
    \]
    with
    \begin{align*}
        D_k &= \begin{cases}
        (\omega^n-\omega)N+(\overline{\omega}^n-\overline{\omega})N^T, & k\le u \\
        (\omega^n-\overline{\omega})N+(\overline{\omega}^n-\omega)N^T, & k>u
        \end{cases} \\
        U &= (\omega^{2u-n}-\omega-\overline{\omega}+1)N + (\overline{\omega}^{2u-n}-1)N^T \\
        L &= (\omega^{2u-n}-1)N + (\overline{\omega}^{2u-n}-\overline{\omega}-\omega+1)N^T.
    \end{align*}
    As before, subtracting $\frac{1}{2(n-2)}$ times each block row and column from every other block row and column will result in a Kronecker product of a matrix $S$ and $(N-N^T)$. Now the first $u-1$ diagonal entries of $S$ are $\frac{1}{2}(-\omega+\overline{\omega})$ and the last $n-u$ are $\frac{1}{2}(\omega-\overline{\omega})$. The entries $d$ above a diagonal have three different forms. If the block was right of $D_i$ and above $D_j$ prior to this last step, then the corresponding entry of $S$ is now
    \[ \begin{cases}
    \frac{n-2-d}{2(n-2)}\overline{\omega}+\frac{n-d}{2(n-2)}\omega-\frac{n-1-d}{n-2}, & i\le u\text{ and } j\le u \\
    \frac{n-1-d}{2(n-2)}\overline{\omega}+\frac{n-1-d}{2(n-2)}\omega-\frac{n-1-d}{n-2}, & i\le u\text{ and } j > u \\
    \frac{n-d}{2(n-2)}\overline{\omega}+\frac{n-2-d}{2(n-2)}\omega-\frac{n-1-d}{n-2}, & i>u\text{ and } j>u.
    \end{cases} \]
    The entries below the diagonal take a similar form; note $S$ is skew-Hermitian. Thus we find that the signature of $B$ is equal to $\sgn((1-\omega^{2u-n})N+(1-\overline{\omega}^{2u-n})N^T)$.
    
    In the case $\varepsilon=-1$, the same steps give that the signature is $\sgn((1-\omega^n)N-(1-\overline{\omega}^n)N^T)+\sgn(S)\sgn(N+N^T)$ as expected, but once again it does not seem there is an easy way to compute $\sgn(S)$. It turns out (from \cite{cha}) that $\sgn(S)$ has the same formula as before, just replacing $n$ with $2u-n$.
\end{itemize}

\bibliographystyle{alpha}
\bibliography{biblio}

\end{document}